\documentclass[final]{siamltex}

\usepackage{amsfonts}
\usepackage{amssymb,amsmath}
\usepackage[mathscr]{eucal}
\usepackage{graphicx}
\usepackage{times}
          
\title{Spectral Shape Preserving Approximation}

\author{Vladimir S. Chelyshkov\thanks{Eastern Kentucky University}}

\begin{document}

\maketitle

\begin{abstract}
We introduce an algorithm of joint approximation of a function and its first derivative by alternative orthogonal polynomials on the interval $[0,1]$. The algorithm exhibits properties of shape preserving approximation for the function.  A weak formulation of approximation is presented. An example on shape preserving extrapolation is given. The weak form is reduced to approximation on a discrete set of abscissas. 

Also, we introduce a new system of orthogonal functions with nice properties -- structured orthogonal polynomials -- and show that the system can be employed for a different kind of joint approximation of a function and its first derivative and may have property of shape preserving approximation. In addition, we show that structured orthogonal polynomials generate wavelet functions.

We complement these results with definition of structured semiorthogonal polynomials and introduce  wavelet basis functions.
\end{abstract}



\pagestyle{myheadings}
\thispagestyle{plain}
\markboth{V.~S.~CHELYSHKOV}{SHAPE PRESERVING APPROXIMATION}

\section{Introduction} Theory  of shape preserving approximation by polynomials has been developing intensively, and we refer here to survey \cite{K} and monograph \cite{G} for most recent assembled results.  In these publications, among topics discussed, much attention was paid to degree of shape preserving approximation with various constrains, and a number of Bernstein-type operators were described (see also paper \cite{K1}). 

Bernstein polynomials \cite{B} approximate  uniformly any continuous function and provide shape preserving approximation. Bernstein operator is linear and monotone. Approximation by the polynomials admits differentiations, and convergence of the polynomial expansion even for smooth functions is slow.

In this paper, a distinct, less general approach for shape preserving approximation is developed. 
We use joint approximation of a function and its first derivative by alternative orthogonal polynomials  \cite{VC1}  (see also  \cite{VC6}) to formally construct an operator for shape preserving approximation  of the function. 
Evidently, the operator is not valid to represent the Weierstrass approximation theorem, it  has $n$-th degree polynomial reproduction property and  exhibits faster convergence for smooth functions.

A few numerical examples on approximation are given, but rigorous
mathematical arguments qualifying our statements on shape preserving approximation and on the rate of convergence are not  presented in this paper.

The paper is self-contained for numerical implementation of the algorithm.

\section{A-kind system} We use a system of alternative orthogonal polynomials
\begin{equation}
\mbox{\boldmath ${\mathscr A}$}_{n}(x)=
\{{\mathscr A}_{nk}(x)\}_{k=n}^{0}
\label{w0}
\end{equation}
for constructing the approximation algorithm. The system has nice properties, which are similar to properties of the classical orthogonal polynomials    \cite{VC1}. Some properties are shown below.    
 
The polynomials ${\mathscr A}_{nk}(x)$ obey the orthogonality relations
\begin{equation}
\int_0^1\frac{1}{x}
{\mathscr A}_{nk}(x)
{\mathscr A}_{nl}(x){\rm d}x=\frac{\delta_{kl}}{k+l},\quad k=n, n-1, ..., 0,\quad l=n, n-1,...,1,
\label{w1}
\end{equation}
but the polynomial ${\mathscr A}_{n0}(x)$ is not normalizable with the given weight. Thus, (\ref{w0}) is a marginal system that contains a singular term ${\mathscr A}_{n0}(x)$. From properties of the system it follows that ${\mathscr A}_{n0}(x)$ are shifted to the interval $[0,1]$ Legendre polynomials. 

The polynomials can be calculated by 
the three-term recurrence relation
\[
{\mathscr A}_{nn}(x)=x^{n},\quad {\mathscr A}_{n,n-1}(x)=(2n-1)x^{n-1}-2nx^{n},
\]
\[
(2k+1)(n+k)(n-k+1){\mathscr A}_{n,k-1}(x)
\]
\[
=2k[(2k-1)(2k+1)x^{-1}-2(n^{2}+k^{2}+n)]{\mathscr A}_{nk}(x)
\]
\[
-(2k-1)(n-k)(n+k+1){\mathscr A}_{n,k+1}(x).
\]
Also, 
\begin{equation}
{\mathscr A}_{nk}^{\prime}(x)=k\frac{{\mathscr A}_{nk}(x)}{x}+2\sum\limits_{l=k+1}^{n}(-1)^{l-k}l\frac{{\mathscr A}_{nl}(x)}{x},\quad k=1, ..., n,
\label{w2}                                                 
\end{equation}
and 
\begin{equation}
{\mathscr A}_{nk}(1)=(-1)^{n-k}.
\label{wplus}
\end{equation}

The system $\mbox{\boldmath ${\mathscr A}$}_{n}(x)$ generates an alternative Gauss-type quadrature with the given weight function. The quadrature is exact for $x^{m}$, $1\le m\le 2n$ on the interval $[0,1]$. This results in the second (discrete) orthogonality property
\begin{equation}
\sum_{j=1}^{n}\frac{w_{j}}{x_{j}}{\mathscr A}_{nk}(x_j){\mathscr A}_{nl}(x_j)=\frac{\delta_{kl}}{k+l},\quad k,l=n, n-1, ..., 1
\label{w3}
\end{equation}
where $w_j$ and $x_j$ are the weights and abscissas of the shifted Legendre-Gauss quadrature.
\begin{figure}[htbp]
\centerline{\includegraphics[height=36mm, width=60mm]{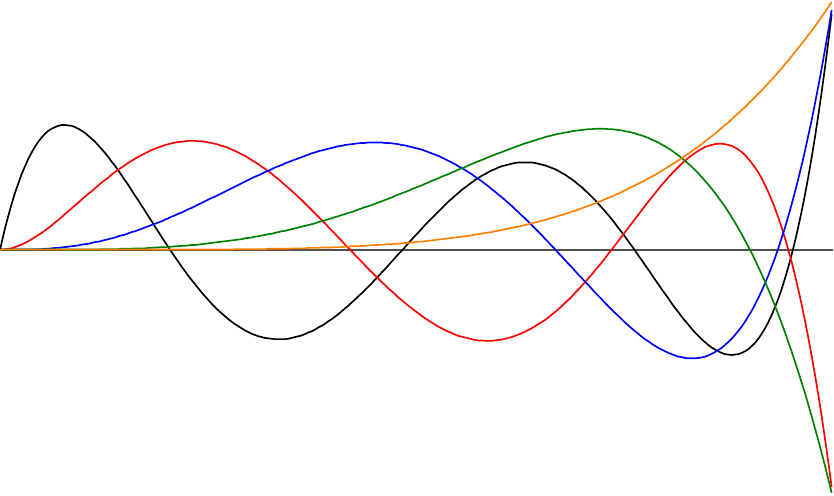}}
\caption{Basis functions: $n=5, k=1-5$.}
\normalsize
\end{figure}

We form an almost orthogonolized polynomial basis 
\begin{equation}
\bold{A}_{n}(x)=\{1, {\mathscr A}_{nk}(x)\}_{k=n}^{1}   
\label{basis}
\end{equation}
in vector space $\Pi_n$ of polynomials of degree less or equal than $n$    and introduce a set
\begin{equation}
\left(\bold{A}_{n} \cup {\mathscr A}_{n0}\right) (x).
\label{2w77}
\end{equation}
for function approximation  on the interval $[0,1]$.

Making use of the basis  $\bold{A}_{n}(x)$ and orthogonality (\ref{w1}), or discrete orthogonality (\ref{w3}) on the abscissas $x_{j}$ of ${\mathscr A}_{n0}(x)$, 
 one can approximate a function by, correspondingly, minimization error in 
a space of quadratically integrable  functions with interpolation at $x=0$, or interpolation  at $n+1$ points $\{0,  x_{j}\}_{j=1}^{n}$.

\section{Spectral approximation}\hspace{-4pt}\footnote{We identify the meaning of the term ``spectral approximation'' with that one  given in \cite{CHQZ},  page 31.}
We initially suppose that $f(x)\in C^{1}([0,1])$ and
expand the function $f(x)$  in sums using the basis ${\bold A}_{n}(x)$ as follows
\[
f(x)=f(0)+f_{0}(x),\quad 
f_{0}(x)=\lim_{n\rightarrow\infty}{f}_{n}(x), 
\]
where
\begin{equation}
f_{n}(x) =\sum_{k=1}^{n}a_{nk}{\mathscr A}_{nk}(x), \quad  a_{nk}\equiv a_{nk}(f_0).
\label{3w1}
\end{equation}
To find coefficients $a_{nk}$ we consider finite projection of
\[ 
f_{0}^{\prime}(x)=\lim_{n\rightarrow\infty}{f}_{n}^{\prime}(x)=\lim_{n\rightarrow\infty}\sum_{k=1}^{n}a_{nk}{\mathscr A}_{nk}^{\prime}(x)
\]
on $\{{\mathscr A}_{nl}(x)\}_{l=1}^{n}$. Making use of  (\ref{3w1}), (\ref{w2}) and (\ref{w1}) we get a system of linear
equations
\begin{equation}
 \sum_{k=1}^{n}T_{lk}a_{nk}=b_{nl},
\label{3w3}
\end{equation}
where
\[
T_{lk}=
\left \{
\begin{array}{c}
\,\,\,\,\,\,\,\,\,\,\,\,\,\,\,\,\,\,0, \\ [1ex]
\,\,\,\,\,\,\,\,\,\,\,1/2,   \\ [1ex]
(-1)^{l-k},
\end{array}
\right.
\begin{array}{c}
l<k, \\ [1ex]
l=k,	\\ [1ex]
l>k
\end{array}
\]
and 
\begin{equation}
b_{nl}=\int_{0}^{1}f_{0}
^{\prime}(x) {\mathscr A}_{nl}(x){\rm d}x, \quad b_{nl}\equiv b_{nl}(f_{0}^{\prime}).
\label{3w77}
\end{equation}
Solving equations (\ref{3w3}) for $a_{nk}$, substituting the solution to (\ref{3w1}), and
applying simple transformations we obtain
\[
f_{n}(x)= 2\sum_{k=1}^{n}b_{nk}{\mathscr B}_{nk}(x),
\]
where
\begin{equation}
{\mathscr B}_{nk}(x)={\mathscr A}_{nk}(x)+2\sum\limits_{l=k+1}^{n}{\mathscr A}_{nl}(x).
\label{3w5}
\end{equation}
Additionally, from (\ref{3w5}) it follows
\begin{equation}
{\mathscr B}_{n0}(x)=1.
\label{3w56}
\end{equation}
It is worthy of notice
\begin{equation}
{\mathscr B}_{nk}^{\prime}(x)=k\frac{{\mathscr A}_{nk}(x)}{x},\quad k=1,2, ..., n,
\label{3w6}
\end{equation}
and we find that (\ref{3w5}) and (\ref{3w56}) form integral co-basis $\bold{B}_{n}(x)=\{1,\, 2 {\mathscr B}_{nk}(x)\}_{k=1}^{n}$ in $\Pi_n$ with
\[
{\mathscr B}_{nk}(x)=k\int_{0}^{x}\frac{{\mathscr A}_{nk}(t)}{t}{\rm d}t,\quad k=1,2,..., n.
\]
Essentially,  basis $\bold{B}_{n}(x)$ is not orthogonal, but the system $\bold{B}_{n}^{\prime}(x)$  is orthogonal
\[
\int_{0}^{1}x{\mathscr B}_{nk}^{\prime}(x){\mathscr B}_{nl}^{\prime}(x){\rm d}x=\frac{k+l}{4}\delta_{kl},\quad 
k,l =1,2,..., n,
\]
and the systems  $\mbox{\boldmath ${\mathscr A}$}_{n}(x)$  and  $\bold{B}_{n}^{\prime}(x)$ have a property of one degree shifted orthogonality with weight function 1
\[
\int_{0}^{1}{\mathscr A}_{nk}(x){\mathscr B}_{nl}^{\prime}(x){\rm d}x=\frac{\delta_{kl}}{2},\quad k,l =0,1,2,..., n.
\]
\begin{figure}[htbp]
\centerline{\includegraphics[height=36mm, width=60mm]{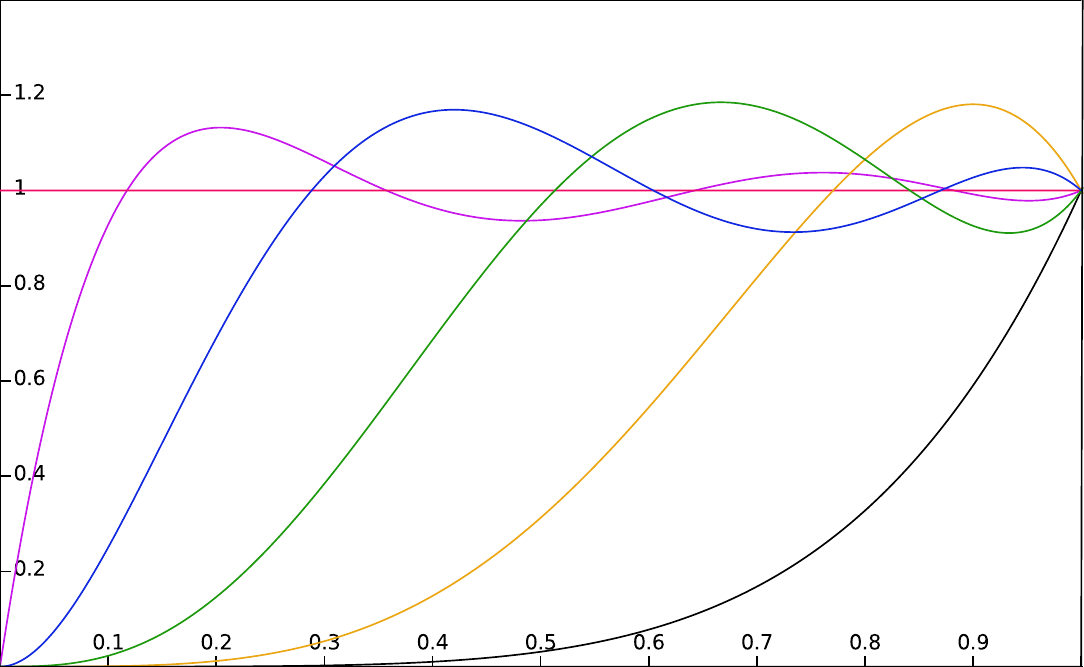}}
\caption{
$\protect{\mathscr B}_{nk}(x),$
 $n=5, k=0-5$.}
\normalsize
\end{figure}

From (\ref{3w6}) and (\ref{w1}) it follows that (\ref{3w7}) is the result of expansion of $f^{\prime}(x)$ in ${\mathscr A}_{nk}(x)/x$, and
\begin{equation}
f(x) =f(0)+ 2 \lim_{n\rightarrow\infty}\sum_{k=1}^{n}b_{nk}{\mathscr B}_{nk}(x),
\label{3w4}
\end{equation}
\begin{equation}
f^{\prime}(x)=2\lim_{n\rightarrow\infty}\sum\limits_{k=1}^{n}b_{nk}{\mathscr B}_{nk}^{\prime}(x).
\label{3w7}
\end{equation}
is joint approximation of the function $f(x)$ and its first derivative by 
the basis $\bold{A}_{n}(x)$ and the integral co-basis $\bold{B}_{n}(x)$ \footnote{At this point, the algorithm is suggestive of a spectral method for solving an initial value problem.}. 

Let us consider three examples on function approximation on the interval $[0,1]$ by expansion (\ref{3w4}). 

For $n=3$ we have
\[
\mbox{ln}(1+x)\approx(342-492\mbox{ln}2)x-(645-930\mbox{ln}2)x^2+(\tfrac{1040}{3}-500\mbox{ln}2)x^3,
\]
\[
1-\mbox{sin}(\pi x)\approx 1+\tfrac{12}{\pi^3}((17\pi^2-180)x-(35\pi^2-360)x^2+(20\pi^2-200)x^3),
\]
for $n=5$
\[
\sqrt{x}\approx\tfrac{2}{11}(15x-35x^2+56x^3-45x^4+14x^5).
\]
The examples represent low degree approximation of a monotonic function, of an even convex function by a polynomial of odd degree, as well as approximation of a monotonic function that is not differentiable at the left end of the interval $[0,1]$. Graphs of all the three expose shape preserving approximation.

One may state that the expansions admit one time differentiation as $n\rightarrow\infty$, and the derivatives of the functions and their approximations have points of intercepts. For, say, $n=3$ and $f(x)=1- \mbox{sin}(\pi x)$  the points are $x_{1}=$ $0.213063$, $x_{2}=0.585763$, and $x_{3}=0.907986$. Since $x_{1}$, $x_{2}$ and $x_3$ cannot be determined  unless the expansion is obtained, distribution of the nodes in the interval $(0,1)$ depends on shape of the derivative at large.

In this paper, we use the term ``concealed interpolation" for approximation that results in interpolation of a function in $O(n)$ non-preassign nodes as $n\rightarrow \infty$.

\section{Weak formulation, composition and asymmetry} Joint approximation\,\,\,  (\ref{3w4})~--~(\ref{3w7}) can be expressed in terms of $f(x)$ as follows. 

We introduce the operator
\begin{equation}
\Omega_{n}(f;x) := f(0)+  2\sum_{k=1}^{n}k\int_{0}^{x}\frac{{\mathscr A}_{nk}(s)}{s}{\rm d}s\int_{0}^{1}f_{0}
^{\prime}(t) {\mathscr A}_{nk}(t){\rm d}t
\label{6w0}
\end{equation}
that corresponds to approximation (\ref{3w4}).  
Integrating second integral in (\ref{6w0}) by parts and applying (\ref{w2})  and (\ref{wplus}) we find
\begin{equation}
\int_{0}^{1}f_{0}
^{\prime}(t) {\mathscr A}_{nk}(t){\rm d}t=(-1)^{n-k}\,f_{0}(1) - c_{nk}/2 - \sum\limits_{l=k+1}^{n}(-1)^{l-k}c_{nl},
\label{6w1}
\end{equation}
where
\begin{equation}
c_{nk}=2k\int_{0}^{1}\frac{1}{t}{f_{0}}(t){\mathscr A}_{nk}(t){\rm d}t, \quad k=1, ..., n.
\label{6w2}
\end{equation}
Making use of (\ref{6w1}) we explicate operator (\ref{6w0}) in the weak form
\[
\Omega_{n}
(f;x) := f(0) + \sum_{k=1}^{n}a_{nk}{\mathscr A}_{nk}(x)
\]
with
\begin{equation}
a_{nk}=(-1)^{n-1}\cdot 2f_{0}(1)+ \sum_{l=1}^{n}S_{kl}c_{nl},
\label{6w22}
\end{equation}
and
\begin{equation}
S_{kl}=
\left \{
\begin{array}{c}
\,\,\,\,\,\,\,\,\,\,\,\,\,\,\,-1, \\ [1ex]
\,\,\,\,\,\,\,\,\,\,\,\,\,\,\,\,\,\,\,3, \\ [1ex]
(-1)^{l}\cdot 2,
\end{array}
\right.
\begin{array}{c}
\,k=l\,\,\, \,\mbox{odd}, \\ [1ex]
\,\,\,k=l\,\,\,\, \mbox{even},\\ [1ex]
 k \ne l.
\end{array}
\label{6w221}
\end{equation}
That is,
\begin{equation}
\Omega_{n}
(f;x):= f(0) +2\sum_{k=1}^{n}{\mathscr A}_{nk}(x)\Big((-1)^{n-1} f_{0}(1)+ \sum_{l=1}^{n}S_{kl}l\int_{0}^{1}\frac{1}{t}{f_{0}}(t){\mathscr A}_{nl}(t){\rm d}t\Big).
\label{6w21a}
\end{equation}
Also, following property of orthogonality (\ref{w1}), we find that evaluation of coefficients $c_{nk}$ results in approximation of $f_{0}(x)$ by
\begin{equation}
\varphi_{n}(x) :=\sum_{k=1}^{n}c_{nk}{\mathscr A}_{nk}(x)
\label{6w20}
\end{equation}
in weighted ``conditional'' $L_2[0,1]$, and  (\ref{6w20}) yields the operator
\begin{equation}
\widehat{\Omega}_{n}
(f;x) := f(0)+  2\sum_{k=1}^{n}{\mathscr A}_{nk}(x)\,k\int_{0}^{1}\frac{1}{t}{f_{0}}(t){\mathscr A}_{nk}(t){\rm d}t.
\label{6w201}
\end{equation}
We can state now that mapping coefficients $c_{nk}$ to $a_{nk}$ by linear transformation (\ref{6w22}),  (\ref{6w221})  forms a composition of approximations  $\widehat{\Omega}_{n}(f;x)$  and ${\Omega}_{n}(f;x)$.

Let us consider an example on the composition. For $n=4$ operator (\ref{6w201}) results in expansion
\begin{equation}
\mbox{sin}(\pi x)\approx \tfrac{6(3\pi^2-28)}{\pi^3}{\mathscr  A}_{41}(x)+\tfrac{4}{\pi}{\mathscr  A}_{42}(x)-\tfrac{6(\pi^2-20)}{\pi^3}{\mathscr A}_{43}(x)+\tfrac{8(\pi^2-6)}{\pi^3}{\mathscr A}_{44}(x),
\label{6w2002} 
\end{equation}
and we observe that coefficients in (\ref{6w2002}) represent eigenvector of matrix  (\ref{6w221}) with eigenvalue 1. In general, this pattern of approximation can be described as follows.

Let $f_{0}(x)$ be an even/odd function, and $n$ is an even/odd number. Then 
\[
\widehat{\varphi}_{n}(x)  := \widehat{\Omega}_{n}(f_{0};x)
\]
is  symmetric concealed interpolation of $f_{0}(x)$ and
\begin{equation}
{\Omega}_{n}(\widehat{\varphi}_{n};x) \equiv \widehat{\varphi}_{n}(x).
\label{6w1017}
\end{equation}
We find that identity (\ref{6w1017}) is reminiscent of our primary choice of approximation of $f(x)$ in the form
\[
f(x )\equiv f(0)+\int_{0}^{x} f^{\prime}(x)\mbox{d}x, \quad f(x) \in C^{1}[0, 1],
\]
and the algorithm does not provide symmetric shape preserving approximation for symmetric $f_{0}(x)$.

Remarkably, $\{{\mathscr A}_{nk}(x)\}_{k=n}^{1}$  is the system of asymmetric functions, and one can change the given pattern of approximation  by distorting  its  symmetry, that is, by choosing odd/even degree of approximation $n$ for even/odd function $f_{0}(x)$. Then  $\widehat{\Omega}_{n}(f_{0};x)$ represents asymmetric concealed interpolation, and ${\Omega}_{n}(f_{0};x)$ stands for asymmetric shape preserving approximation. Hereinafter we consider this option for $n$  as a component of the algorithm.
Accordingly, no special choice of $n$ is required for an asymmetric $f_{0}(x)$, for any asymmetric function can be represented as a sum of even and odd function.

With the pattern selected, we find that the composition of approximations consists of two different types of approximation by the same system of functions.

It was mentioned above that expansion of $f(x)=\sqrt{x}$ by operator  ${\Omega}_{n}(f;x)$ is one time differentiable, and  it is easy to verify that expansion of the same function by $\widehat{\Omega}_{n}(f;x)$ is not. 

We infer  that operator (\ref{6w21a}) may  be  applied for approximation of functions from a wider class, say, $f(x)\in C\cap BV([0,1])$, as well as for approximation of the derivative $f^{\prime}(x)$ where it exists on $[0,1]$ \footnote{Below we suppose that the second integral in (\ref{6w0}) 
is Henstock - Kurzweil integral.}.

\section{An example on shape preserving extrapolation}In this section we numerically compare the two 
types of approximation of a particular function of $C^{\infty}([0,1]$) class. We choose
\[
f(x)=\mbox{sin}(\tfrac{\pi}{2}x),\quad f_{1}(x) = \widehat{\Omega}_{9}
(f;x),\quad f_{2}(x) = \Omega_{9}
(f;x), \quad x \in [0, 1].
\]
Numerical results for $n=9$ are presented in Fig. 7.1. The function $f_1(x)$ (yellow curve) represents $L_{2}$-approximation with concealed nodes of interpolation in the interval $(0,1)$ and the node at $x=0$, whereas shape preserving approximation $f_2(x)$ (green curve) intersects $f(x)$ only at $x=0$. The two curves are visually very close to the graph of $\mbox{sin}(\tfrac{\pi}{2}x)$ (red curve) in a much longer interval, say, $[-1, 2]$. 
 One may also observe that shape preserving approximation provides a ``better'' extrapolation, for it is in a compliance with the $f(x)$ if $x \in [-1.5, -1 ]\cup [2,2.5]$.  We see that shape preserving approximation may lead to a better result near turning points outside of the original interval.

The observation is more of general interest, rather than of computational implementation, 
\begin{figure}[h!]
\centerline{\includegraphics[height=36mm, width=60mm]{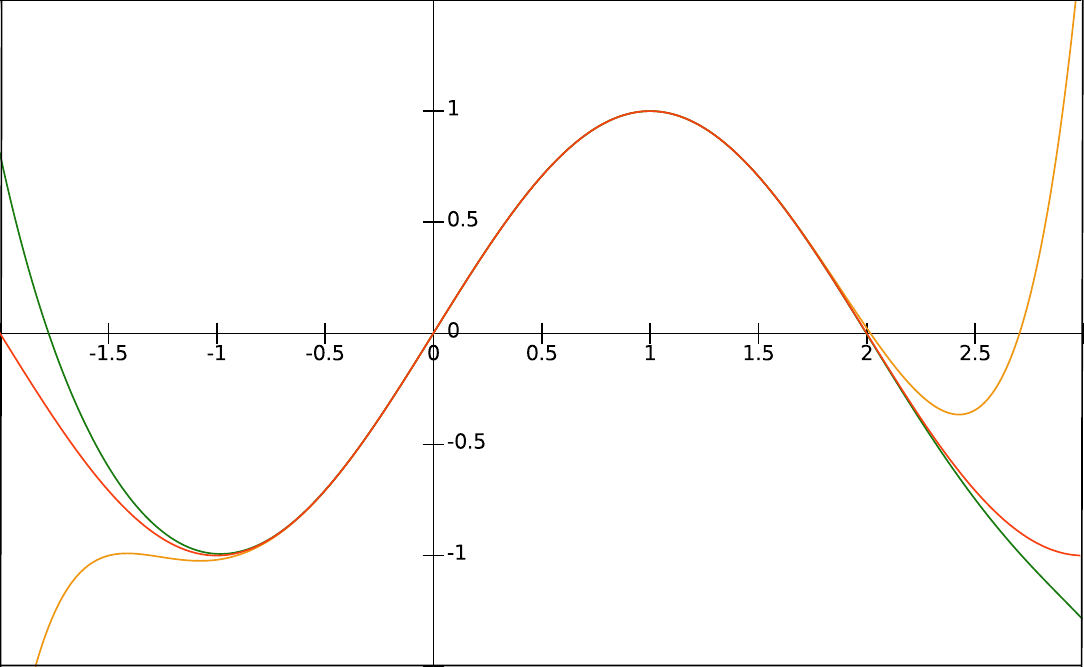}}
\caption{Extrapolation by the composition.}
\normalsize
\end{figure} 
since calculation of orthogonal polynomials outside of the interval of orthogonality meets difficulty related to a number of significant figures required. 

\section{Discretization and pseudo-basis representation} Let us consider a discrete weak formulation of joint approximation. Lagrange interpolation  in the Gaussian nodes of a function from the above selected class converges uniformly \cite{PV}, and we employ discrete orthogonality property (\ref{w3}) to interpolate the $f_{0}(x)$ by the polynomial 
\begin{equation}
\label{6w3}
{\psi}_{n}(x):=\sum_{k=1}^{n}d_{nk}{\mathscr A}_{nk}(x)
\end{equation}
with
\[
d_{nk}=2k\sum_{i=1}^{n}\frac{w_{j}}{x_{j}}f_{0}(x_{j}){\mathscr A}_{nk}(x_j).
\]
Follwing (\ref{6w3}) we introduce the operator
\begin{equation}
\widehat{\mathscr W}_{n}
(f;x) :=f(0)+ 2\sum_{j=1}^{n}f_{0}(x_{j}){\mathscr Q}_{nj}(x),
\label{6w32}
\end{equation}
\vspace{-5mm}
\[
{\mathscr Q}_{nj}(x) :=\frac{w_{j}}{x_{j}}\sum_{k=1}^{n} k {\mathscr A}_{nk}(x_j){\mathscr A}_{nk}(x)
\]
that interpolates $f(x)$ in $n$ Gaussian nodes and in $x=0$.

Making use of the interpolant ${\psi}_{n}(x)$ for evaluating integrals 
in (\ref{6w21a}) we finally obtain  the operator of discrete approximation of $f(x)$
\[
\hspace{-40mm}{\mathscr W}_{n}(f;x):=f(0)+2\sum_{k=1}^{n}{\mathscr A}_{nk}(x)\Big((-1)^{n-1} f_{0}(1)
\]
\vspace{-8pt}
\begin{equation}
\hspace{42mm}+\sum_{l=1}^{n}S_{kl}l\sum_{j=1}^{n}\frac{w_{j}}{x_{j}}f_{0}(x_{j}){\mathscr A}_{nl}(x_j)\Big).
\label{6w4}
\end{equation}

Again, linear transformation (\ref{6w22}), being applied to $d_{nk}$,
maps interpolation  (\ref{6w32}) to shape preserving approximation (\ref{6w4}), the operators form a composition of two discrete approximations by the same system of polynomials, and the nodes of interpolation are abscissas for the joint approximation.

Formalism (\ref{6w21a}), (\ref{6w3}), ({\ref{6w4})   represents exact projection that unites weak formulation and collocation method; it was introduced as compound spectral formalism in \cite{VC2} for solving a problem on linear theory of hydrodynamic stability in ordinary differential equations.

Operator (\ref{6w4}) can be expressed in similar to (\ref{6w32}) form  as
\begin{equation}
{\mathscr W}_{n}
(f;x)=f(0){\mathscr P}_{n0}(x)+\sum_{j=1}^{n}f_{0}(x_{j}){\mathscr P}_{nj}(x)+f_{0}(1){\mathscr P}_{n,n+1}(x),
\label{6ww1}
\end{equation}
\vspace{-10pt}
\[
{\mathscr P}_{n0}(x)=1, \quad {\mathscr P}_{n,n+1}(x)=2(-1)^{n-1}\sum_{k=1}^{n}{\mathscr A}_{nk}(x),
\]
\vspace{-10pt}
\[
{\mathscr P}_{nj}(x)=2\frac{w_{j}}{x_{j}}\sum_{k=1}^{n}{\mathscr A}_{nk}(x) \sum_{l=1}^{n}S_{kl}l{\mathscr A}_{nl}(x_j).
\]
where ${\{\mathscr P}_{nj}(x)\}_{j=0}^{n+1}$ is a pseudo-basis, and ${\mathscr P}_{n,n+1}(x)$ is considered as a linear combination of ${\mathscr P}_{nj}(x)$, $0\le j \le n$. Approximation  by (\ref{6ww1}) possesses $n$-th degree polynomial reproduction property; it interpolates  $f(x)$ only at $x_{0}=0$, and approximation of the first derivative of $f(x)$ is a concealed interpolation.  

The polynomials ${\mathscr P}_{nj}(x)$ $j=0,1,...,n$ are not orthogonal. These polynomials have no zeros in the interval $(0, 1]$, increase in absolute value for $j>0$ as $n$ increases, and alternate in sign with respect to $j$ in $[0,1]$. That is, representation of ${\mathscr W}_{n}(f;x)$ in form of (\ref{6ww1}) may result in numerical instability for greater $n$, and the algorithm of discrete approximation should be implemented in its original form (\ref{6w4}).

\section{Convergence}
Obviously, approximation of $f_{0}(x)$  by $\varphi_{n}(x)$  may not be distinguished as a separate part of the algorithm; furthermore, interpolation of $f_{0}(x)$ by $\psi_{n}(x)$ is not the necessary step of the discrete variant of the algorithm. Indeed, coefficients $d_{nk}$ are the result of evaluation of $c_{nk}$ by the $n$-th order alternative Gauss quadrature in (\ref{6w2}). So, both continuous and discrete variants of the algorithm may be considered  without introducing the compositions. 

At this point we recognize the algorithm as a whole entity and assume that the most general class of functions for approximation of $f(x)$ might be defined in the following way. 

$Conjecture.$  Let $f(x)$ be a continuous almost everywhere differentiable function in the interval $[0, 1]$. Then approximation of $f(x)$ by the operator $\Omega_{n}(f;x)$ in weak formulation (\ref{6w21a})  and  its discrete analog ${\mathscr W}_{n}(f;x)$ converges uniformly in $[0, 1]$.

\section{Conclusion} Presented algorithm of joint approximation does not provide interpolation of a function at the right end of the interval, but example on extrapolation of a smooth function shows that low degree approximation may be good enough in the closed interval $[0,1]$ without interpolation at the endpoint.  
Thus, the algorithm may originate polynomial curves for solving problems of isogeometric analysis \cite{Iso}.

Formally, the algorithm can be reconstructed for shape preserving approximation on a half-line, and it can be done in two different ways. 

Firstly, the system of orthogonal exponential functions  ${\mathscr E}_{nk}(t):={\mathscr A}_{nk}(e^{-t})$, $n>0$, described in \cite{VC1, VC6}  can be employed for reformulation of the algorithm. One can find that the rearranged version of the algorithm provides  approximation for functions $f(t) \in C^{1}[0, \infty)$,  $f(t) \sim e^{-\alpha t}$ as $t \rightarrow \infty$,\, and $\alpha>0$. 

Secondly,  reformulation can be performed by choosing  introduced in \cite{VC7} system of orthogonal rational functions ${\cal R}_{nk}^{{\mathscr A}}(t):={\mathscr A}_{nk}(t^{-1})$. Then, this version of the algorithm can be employed for approximation of functions $f(t) \in C^{1}[1, \infty)$ and $f(t) \sim t^{-\alpha}$ as $t \rightarrow \infty$.

\section{Appendix 1. Structured orthogonal polynomials and an example on approximation. Preliminaries}Standard theory of orthogonal polynomials on an interval originates from Gram-Schmidt process of orthogonalization of the sequence of monomials in the order of the exponent increase.  It results in sequences of polynomials that  have distinct interlaced real zeros in an  interval $(a,b)$
and other fundamental properties  \cite{CHI, IM}. Alternative orthogonal polynomials \cite{VC1, VCETNA, VC6} are systems of functions introduced by orthogonalization of the sequence of monomials on the interval $[0,1]$ in the inverse order, from $k=n$ to $k=0$. In addition to zeros inside the interval, the  polynomials have $k$-multiple zeros at the left end of the interval, and their other properties are similar to those ones of the standard theory. 

Below we apply alternative orthogonalization  to a specially chosen polynomial sequence of functions $\{\pi_{k}(x)\}_{k=n}^{0}$ and introduce an orthogonal polynomial structure that have $k$ multiple zeros distributed both at the left and right ends of the interval  $[0,1]$.

Specifically, we  orthogonalize the sequence
\begin{equation}
\quad \pi_{k}(x):=x^{k-\left \lfloor{\frac{k}{2}}\right \rfloor} (1-x )^{\left \lfloor{\frac{k}{2}}\right \rfloor},\quad k=0,1, ..., n,\quad n=1,2,...
\label{ap1}
\end{equation}
beginning with $k=n$ in the inverse order until $k=0$ with  weight function $1/x$.
This results in a system of functions
\[
\mbox{\boldmath ${\mathscr S}$}_{n}(x)=
\{{\mathscr S}_{nk}(x)\}_{k=n}^{0},
\]
such that 
\begin{equation}
\int_{0}^{1}\frac{1}{x}{\mathscr S}_{nk}(x){\mathscr S}_{nl}(x){\rm d}x=\delta_{kl},\quad k=n,n-1,...,0,\quad l=n,n-1,...,1
\label{ap2}
\end{equation}
Obviously, polynomial ${\mathscr S}_{n0}(x)$ cannot be normalized by (\ref{ap2}), but it is orthogonal to ${\mathscr S}_{nk}(x)$ for $k>0$.
We also may interpret definition of orthogonality relations (\ref{ap2}) as one degree shifted orthogonality with weight function 1. That is,
\[
\int_{0}^{1}{\mathscr S}_{nk}(x)\hspace{1pt}^{\prime}\hspace{-1pt}{\mathscr S}_{nl}(x){\rm d}x=\delta_{kl}, \quad ^{\prime}\hspace{-1pt}{\mathscr S}_{nl}(x)=\frac{{\mathscr S}_{nl}(x)}{x} .
\]
Remarkably, one can find that structured orthogonal polynomials ${\mathscr S}_{nk}(x)$ farther can be defined by Rodrigues' type formula 
\[
{\mathscr S}_{nk}(x)=
 c_{nk}\frac{\pi_{k}(x)}{x^k(1-x)^k}
\frac{{\rm d}^{n-k}}{{\rm d}x^{n-k}}(x^{n}(1-x)^{n}),
\]
where $c_{nk}$ is a normalization factor. The formula provides an easy way for characterizing other properties of the polynomials, but their full description deserves a separate consideration, and we do not present it here.

Curiously, polynomials ${\mathscr S}_{n0}(x)$ are the shifted Legendre polynomials $\tilde{P}_{n}(x)$,  i.e.,  ${\mathscr S}_{n0}(x)$  generate the Gauss quadrature, and the zeros of ${\mathscr S}_{n1}(x)$ in the interval $(0,1)$ represent shifted Lobatto abscissas. For $n=3$, the orthogonalization process results in
\[
\tilde{P}_{3}(x)=\pi_{0}(x)-2\pi_{1}(x)-10\pi_{2}(x)+20\pi_{3}(x),
\]
and we may suppose  that those standard orthogonal polynomials that are equipped with Rodrigues formula can be effectively calculated for greater $n$ by their representation in the structured form of linear combinations of $\pi_{k}(x)$.

\begin{figure}[htbp]
\centerline{\includegraphics[height=36mm, width=60mm]{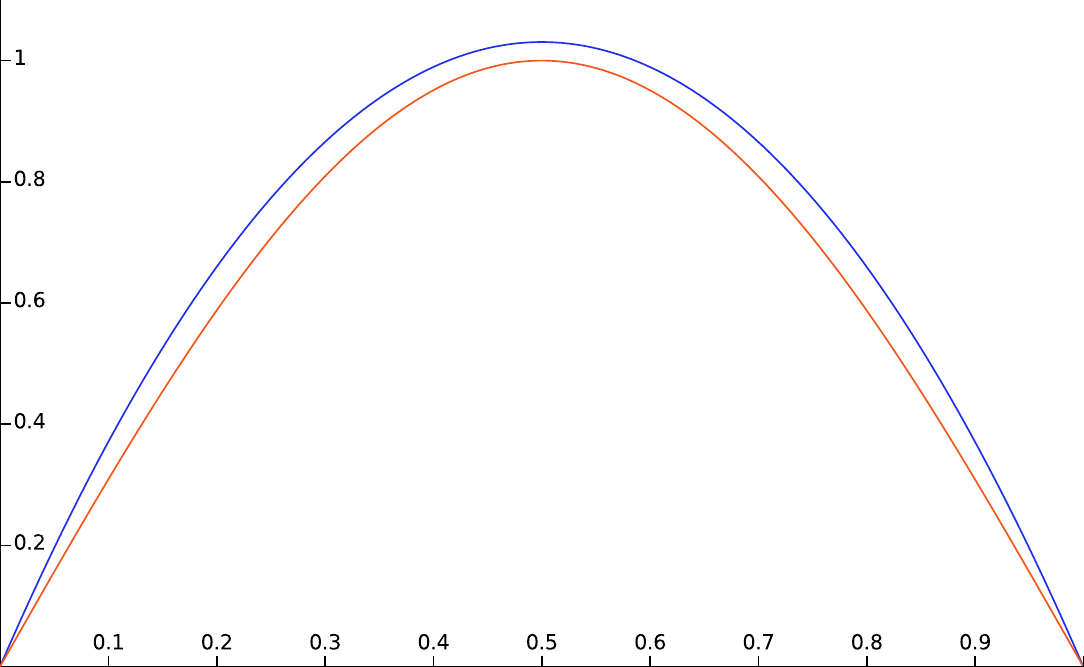}}
\caption{ Function $\sin(\pi x)$ (in red) and its approximation (in blue) for
$n=3$.}
\normalsize
\end{figure}

We believe that ${\mathscr S}_{nk}(x)$ is a good choice for developing next algorithm of joint approximation of a function and its derivative. Making use of this approach for $f(x)=\sin(\pi x)$  on the interval $[0,1]$ with $n=3$ we found 
\[
\sin(\pi x)\approx \tfrac{60(12-\pi^2)}{\pi^3}x(1-x),
\]
that is, the second (sic!) degree shape preserving  polynomial approximation with interpolation in two endpoints (See Fig. 9.1).

\section{Appendix 2. A basis in the space of compactly supported smooth functions and wavelet functions}

Let $x \in \mathbb{R}$, $k>1$, $n\ge k$, $\Lambda_{nk}(x)\equiv{\mathscr S}_{nk}(x)$, if $x \in (0,1)$, and $\Lambda_{nk}(x)\equiv 0$ otherwise. Then system of functions $\Lambda_{nk}(x)$ represents a basis for $L^2(\mathbb{R})$ which consists of compactly supported  functions of class $C^{k-2}$.  The basis is orthonormal with weight function $1/x$ in $\mathbb{R}$.
\begin{figure}[htbp]
\centerline{\includegraphics[height=36mm, width=60mm]{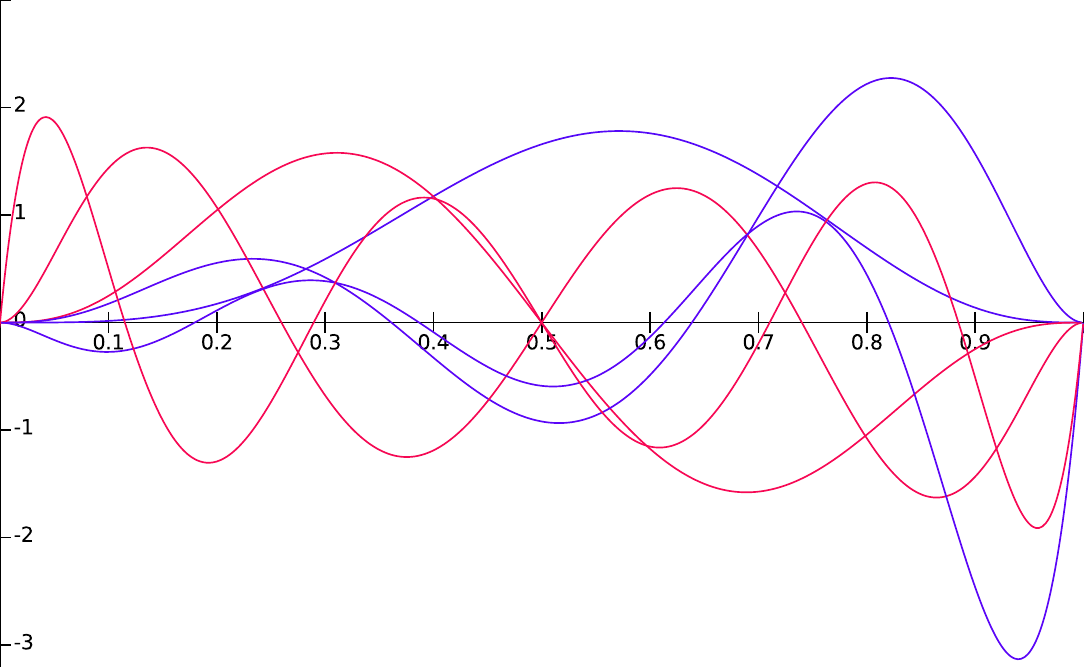}}
\caption{
$\protect{\Lambda}_{nk}(x),$
 $n=7, k=2-7$.}
\normalsize
\end{figure}
Graphs of  $\Lambda_{7k}(x), x \in (0,1), 1<k\le7$,  are  shown  in  Fig 10.1. One may notice that the zeros of $\Lambda_{nk}(x)$ in $(0,1)$ are placed symmetrically, and three of six curves presented  are symmetric with respect to the middle point of the interval (``anti-symmetric''). We find that the set of chosen six functions  is an asymmetric-antisymmetric collection, and  this results in conclusion that antisymmetric subset of  the system $\Lambda_{nk}(x)$ ($n$ is odd, $k $ is even) generates continuous wavelet functions (mother wavelets). Similarly, the case ``$n$ is even, $k$ is odd'' is a proper choice to consider. It should be noted, however, that the system $\Lambda_{nk}(x)$ does not contain symmetric with respect to the line $x=1/2$ sub-sequences that satisfy zero mean condition.

In  Fig. 10.2 we show only antisymmetric part of the system $\Lambda_{nk}(x)$ for $x \in (0,1)$.
\begin{figure}[htbp]
\centerline{\includegraphics[height=36mm, width=60mm]{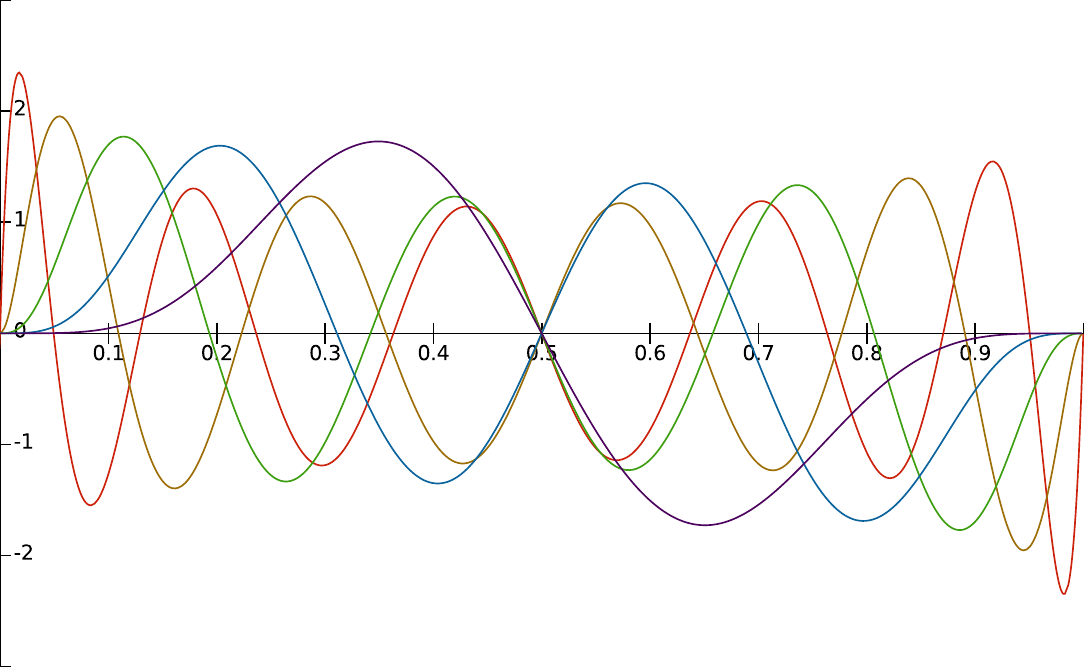}}
\caption{
Wavelets $\protect{\Lambda}_{nk}(x),$
 $n=11, k=2,4,6,8,10$.}
\normalsize
\end{figure}

One may mention that properties of orthogonal polynomials ${\mathscr S}_{nk}(x)$ allow to easily calculate functions $\Lambda_{nk}(x)$ of desirable  smoothness and degree.

\section{Appendix 3. Structured  semiorthogonal polynomials. Preliminaries} 
By giving up orthogonality of ${\mathscr S}_{nk}(x)$ one can introduce a symmetric-antisymmetric system of polynomial functions $\widetilde{{\mathscr S}}_{nk}(x)$ in $L^2(0,1)$.

Let $n$ be odd, and $k=n, n-1, ..., 0$. We define $\widetilde{{\mathscr S}}_{nk}(x)$  in the following way 
\[
\widetilde{{\mathscr S}}_{nk}(x)= d_{nk}\cdot 
\left \{
\begin{array}{c}
\hspace{3pt}{\mathscr S}_{nk}(x)/x,\\  [1ex]
\hspace{-7pt}{\mathscr S}_{nk}(x), \\ [1ex]
\hspace{6pt}\widetilde{\mathscr S}_{n,k+1}^{\prime}(x),  
\end{array}
\right.
\begin{array}{c}
\hspace{-4pt}k=n,\\ [1ex]
k\,\,\,\mbox{even},\\ [1ex]
\hspace{28pt}k\,\,\,\mbox{odd},\, k < n.
\end{array}
\]
For even $n$ we choose 
\[
\widetilde{{\mathscr S}}_{nk}(x)= d_{nk}\cdot 
\left \{
\begin{array}{c}
\hspace{-5pt}{\mathscr S}_{nk}(x),\\  [1ex]
\hspace{6pt}{\mathscr S}_{nk}(x)/x, \\ [1ex]
\hspace{8pt}\widetilde{\mathscr S}_{n,k+1}^{\prime}(x),  
\end{array}
\right.
\begin{array}{c}
\hspace{-4pt}k=n,\\ [1ex]
\hspace{-3pt}k\,\,\,\mbox{odd},\\ [1ex] 
\hspace{31pt}k\,\,\,\mbox{even},\, k < n.
\end{array}
\]
Also, we  select weight function $w(x)=1$ for $\widetilde{{\mathscr S}}_{nk}(x)$ normalization and  define coefficients $d_{nk}$ by condition  
\[
\int_0^1 (\widetilde{{\mathscr S}}_{nk}(x))^2 {\rm d}x=1.
\] 
The system $\{\widetilde{{\mathscr S}}_{nk}(x)\}_{k=n}^{0}$  is a complementary semiorthogonal polynomial structure.                        
Polynomials $\widetilde{{\mathscr S}}_{nn}(x)$ are symmetric  bump functions, or window functions, and, for $0\le k<n$, $\widetilde{{\mathscr S}}_{nk}(x)$ are, correspondingly, odd and even functions in $[0,1]$.

Let $m=2$ if $n$ is odd, and $m=3$ if $n$ is even. For $m\le k\le n$, we introduce system of functions $\widetilde{\Lambda}_{nk}(x)$ in the following way: $\widetilde{\Lambda}_{nk}(x)\equiv \widetilde{{\mathscr S}}_{nk}(x)$,  if  $x\in (0,1)$, and $\widetilde{\Lambda}_{nk}(x)\equiv 0$ otherwise. For $m\le k< n$ functions  $\widetilde{\Lambda}_{nk}(x)$ satisfy the condition for zero mean, i.e., these are wavelet functions.

A compactly supported  function $f(x)\in L^2(0,1)$  may be approximated by the system $\widetilde{\Lambda}_{nk}(x)$.
Let
\begin{equation}
f(x)=\lim_{n\rightarrow\infty}{f}_{n}(x), 
\label{11.1}
\end{equation}
We represent ${f}_{n}(x)$ in the form
\begin{equation}
f_{n}(x) =\sum_{k=n}^{m}a_{nk}\widetilde{\Lambda}_{nk}(x),
\label{11.2}
\end{equation}
reduce the problem (\ref{11.1}), (\ref{11.2}) to solving a system of linear equations for $a_{nl}$
\begin{figure}[htbp]
\centerline{\includegraphics[height=36mm, width=60mm]{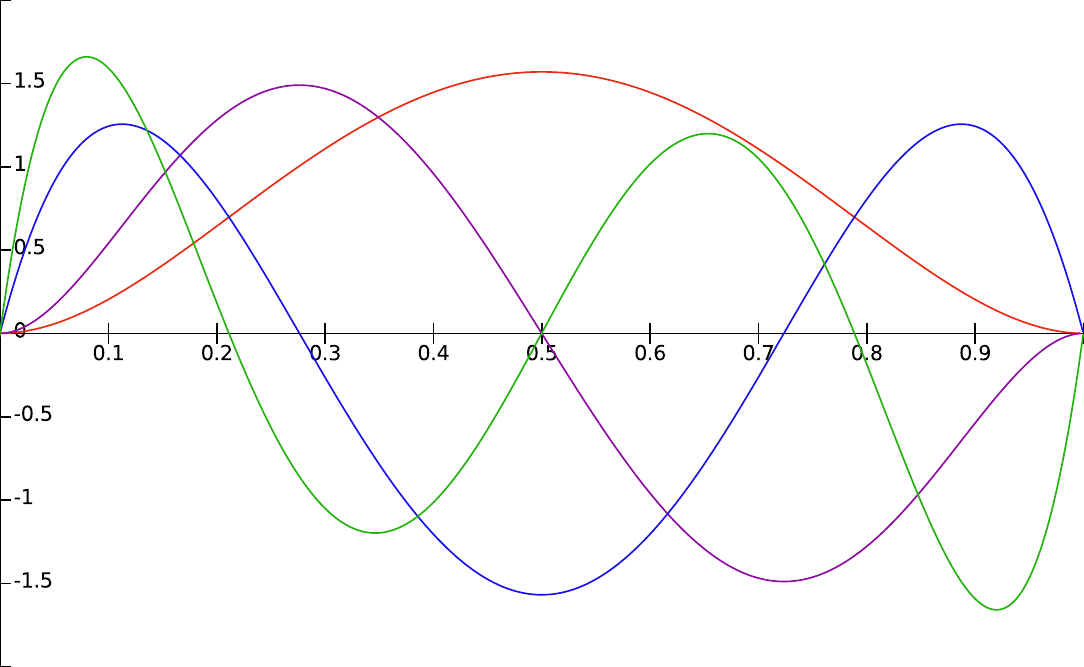}}
\caption{
Functions $\protect{\widetilde{\Lambda}}_{nk}(x),$
 $n=5, k=5-2$.}
\normalsize
\end{figure}
\[
\sum_{l=n}^{m}A_{nkl}a_{nl}=b_{nk},
\]
\[
A_{nkl}=\int_{0}^{1}\widetilde{\Lambda}_{nk}(x)\widetilde{\Lambda}_{nl}(x){\rm d}x,\quad b_{nk}=\int_{0}^{1}f(x)\widetilde{\Lambda}_{nk}(x){\rm d}x,  
\]
and find 
\[
a_{nk}=A_{nkl}^{-1}b_{nl}.
\]

Symmetric matrices  ${\bold A}_{n}={A}_{nkl}$, $k,l=n, ..., m$  have zero-alternate entries due to the functions $\widetilde{\Lambda}_{nk}(x)$ symmetry; also, calculations show that ${\bold A}_{n}$'s parallel to the main diagonal skew rows that contain non-zero entries alternate in sign, and the matrices are positive-definite.  This results in a nice feature
 -- symmetric inverse matrices  ${\bold A}_{n}^{-1}={A}_{nkl}^{-1}$  are non-negative, zero-alternate, and positive-definite. 

Let us consider two examples. For $n=5$ the system 
\[
\widetilde{\Lambda}_{55}(x)=\sqrt{630}x^2(1-x)^2,
\]
\[
\widetilde{\Lambda}_{54}(x)=\sqrt{6930}x^2(1-x)^2(1-2x),
\]
\[
\widetilde{\Lambda}_{53}(x)=\sqrt{630}x(1-x)(1-5x+5x^2),
\]
\[
\widetilde{\Lambda}_{52}(x)=\sqrt{2310}x(1-x)(1-2x)(1-6x+6x^2). 
\]
generates matrices
\[
{\bold A}_{5}=
\begin{bmatrix}
1 & 0 & -\frac{1}{2} & 0  \\
0 & 1 & 0 & -\frac{\sqrt{3}}{6}  \\
-\frac{1}{2} & 0 & 1 & 0 \\
0 & -\frac{\sqrt{3}}{6} & 0 & 1 
\end{bmatrix}
, 
\]
and
\[
{\bold A}_{5}^{-1}=\frac{2}{33}
\begin{bmatrix}
22 &       0       & 11 &       0        \\
 0  &     18       &  0  & 3\sqrt{3} \\
11 &       0       & 22 &       0        \\
 0  & 3\sqrt{3}&  0  &      18 
\end{bmatrix}
. 
\] 
For $n=6$
\[
\widetilde{\Lambda}_{66}(x)=\sqrt{12012}x^3(1-x)^3,
\]
\[
\widetilde{\Lambda}_{65}(x)=\sqrt{6930}x^2(1-x)^2(1-2x),
\]
\[
\widetilde{\Lambda}_{64}(x)=\sqrt{630}x(1-x)(1-5x+5x^2),
\]
\[
\widetilde{\Lambda}_{63}(x)=\sqrt{630}x(1-x)(1-2x)(2-11x+11x^2). 
\]
Thus,
\[
{\bold A}_{6}=
\begin{bmatrix}
1                                   & 0 & -\sqrt{\frac{39}{110}} & 0  \\
0                                   & 1 &                                     0& 0  \\
-\sqrt{\frac{39}{110}}& 0 &                                     1& 0  \\
0                                   & 0 &                                     0& 1 
\end{bmatrix}
, 
\]
and
\[
{\bold A}_{6}^{-1}=\frac{1}{71}
\begin{bmatrix}
110              & 0&\sqrt{4290}&0 \\
 0                 &71&                0&0 \\
\sqrt{4290} &  0&            110&0 \\
 0                 &  0&                0&71
\end{bmatrix} 
.
\]
For $n= 5,6$ matrices ${\bold A}_{n}$ and ${\bold A}_{n}^{-1}$ have positive eigenvalues and orthogonal eigenvectors.  

These and other examples considered  speak in favor of  convergence of expansion (\ref{11.1}), (\ref{11.2}) in $L^2(0,1)$. 
\vspace{4pt}

$Conjecture.$   System of functions $\widetilde{\Lambda}_{nk}(x)$ is a semi-orthogonal (alternative) Riesz basis for compactly supported functions in $L^{2}(0,1)$.
\vspace{6pt}

Surprisingly, the last entry in the second row of  matrix ${\bold A}_6$ is zero. This means that the product of two odd wavelets              $\widetilde{\Lambda}_{65}(x)$ and $\widetilde{\Lambda}_{63}(x)$ forms an even wavelet
\[
_c\widetilde {\Lambda}_{6,5,3}(x)=(14549535/2)^{1/2}x^3(1-x)^3(1-2x)^2(2-11x+11x^2)
\]
with a double zero at $x=1/2$. The wavelet consists of two flipped asymmetric wavelets.

Also, for $n=7$, the fifth entry in the first row of matrix ${\bold A}_7$ is equal to zero, and product of two even wavelets, 
$\widetilde{\Lambda}_{77}(x)$ and $\widetilde{\Lambda}_{73}(x)$,
 forms an even wavelet 
\begin{figure}[htbp]
\centerline{\includegraphics[height=36mm, width=60mm]{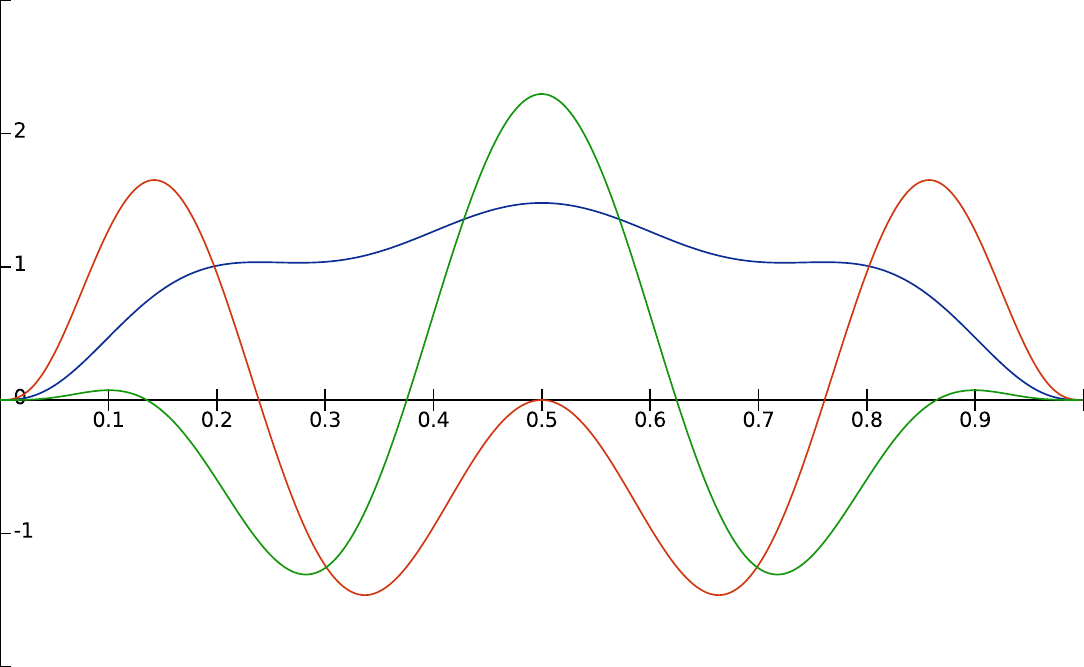}}
\caption{
$\protect{_c\widetilde{\Lambda}}_{6,5,3}(x)$(in red), $\protect{_c\widetilde{\Lambda}}_{7,7,3}(x)  (in\hspace{2pt}green)$, $\protect{_c\widetilde{\Lambda}}_{6,4,3}^{a}(x) (in\hspace{2pt}blue) $}
\normalsize
\end{figure}
\[
_c\widetilde {\Lambda}_{7,7,3}(x)=
(3346393050/1363)^{1/2}x^4(1-x)^4(5-64x+246x^2-364x^3+182x^4).
\]

 The functions $_c\widetilde {\Lambda}_{6,5,3}(x)$ and $_c\widetilde {\Lambda}_{7,7,3}(x)$, being consecutively integrated and differentiated, result in collateral  semi-orthogonal 
sequences of five and six terms respectively.

Similarly, zero entries of  the matrix ${\bold A}_n$ that are the result of an even-odd wavelet product conform to collateral sequences. 

For $n=6$, the fourth entry in the third row of  matrix ${\bold A}_6$ is equal to zero, and the product of two wavelets, 
$\widetilde{\Lambda}_{64}(x)$ and $\widetilde{\Lambda}_{63}(x)$, being integrated, forms a window function
\[
_c\widetilde{\Lambda}_{6,4,3}^{a}(x)=(9699690/271)^{1/2}x^3(1-x)^3(8-63x+195x^2-264x^3+132x^4).
\] 

Collateral sequences  are atypical wavelets and window functions that can be represented by linear combinations of
$\widetilde {\Lambda}_{nk}(x)$ with a properly chosen $n$. 

The system of functions $\widetilde {\Lambda}_{nk}(x)$  is an unbounded collection of wavelets.  It may be of interest for applications in signal analysis \cite{Db}.

\end{document}